\documentclass[12pt]{amsart}
\usepackage{amsmath,amsfonts,amssymb,amsthm}
\usepackage{graphics}
\usepackage{epsfig,psfrag}
\usepackage{comment}

\newtheoremstyle{thm}
  {9pt}{9pt}{\itshape}{}{\bfseries}{}{.5em}{}
\theoremstyle{thm}
\newtheorem{thm}{Theorem}
\newtheorem{cor}[thm]{Corollary}
\newtheorem{lemma}[thm]{Lemma}

\newtheoremstyle{defin}
  {9pt}{9pt}{}{}{\bfseries}{}{.5em}{}
\theoremstyle{defin}

\newtheoremstyle{exm}
  {9pt}{9pt}{}{}{\scshape}{}{.5em}{}
\theoremstyle{exm}
\newtheorem*{exm}{Example}
\newtheorem*{rmk}{Remark}
\newtheoremstyle{proof}
  {}{}{}{}{\itshape}{:}{.5em}{}
\theoremstyle{proof}

\newtheorem*{skt}{Sketch of proof}

\newcommand{\set}[1]{\{#1\}}
\newcommand{\Z}{{\mathbb Z}}
\newcommand{\p}[1]{\mathcal{#1}}
\newcommand{\s}[1]{\mathbf{#1}}

\DeclareMathOperator{\syt}{SYT}
\DeclareMathOperator{\var}{var}

\author{Matja\v z Konvalinka}
\title[Weighted hook-length formula II]{The weighted hook-length formula II: \\Complementary formulas}

\begin{document}

\begin{abstract}
 Recently, a new weighted generalization of the branching rule for the hook lengths, equivalent to the hook formula, was proved. In this paper, we generalize the complementary branching rule, which can be used to prove Burnside's formula. We present three different proofs: bijective, via weighted hook walks, and via the ordinary weighted branching rule.
\end{abstract}

\maketitle
\section{Introduction} \label{intro}

The classical hook-length formula gives an elegant product formula for the number of standard Young tableau. Since its discovery by Frame, Robinson and Thrall in~\cite{FRT}, it has been reproved, generalized and extended in several different ways, and applications 
have been found in a number of fields of mathematics.

\medskip

\begin{figure}[hbt]
\psfrag{x1}{$x_1$}
\psfrag{x2}{$x_2$}
\psfrag{x3}{$x_3$}
\psfrag{x4}{$x_4$}
\psfrag{x5}{$x_5$}
\psfrag{y1}{$y_1$}
\psfrag{y2}{$y_2$}
\psfrag{y3}{$y_3$}
\psfrag{y4}{$y_4$}
\psfrag{y5}{$y_5$}
\psfrag{y6}{$y_6$}
\psfrag{1}{$1$}
\psfrag{2}{$2$}
\psfrag{3}{$3$}
\psfrag{4}{$4$}
\psfrag{5}{$5$}
\psfrag{6}{$6$}
\psfrag{7}{$7$}
\psfrag{L}{$\lambda$}
\begin{center}
\epsfig{file=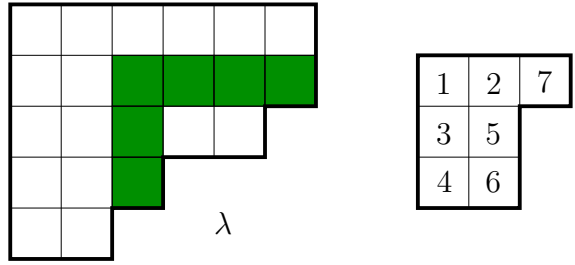,height=3.5cm}
\end{center}
\caption{Young diagram $[\lambda]$,  $\lambda = 66532$, and a hook $H_{23}$ with hook length $h_{23}=6$; a standard Young tableau of shape $322$.}
\label{fig9}
\end{figure}

Let $\lambda = (\lambda_1,\lambda_2,\ldots,\lambda_\ell)$, $\lambda_1 \geq \lambda_2 \geq \ldots \geq \lambda_\ell > 0$, be a partition of $n$, $\lambda \vdash n$, and let $[\lambda] = \set{(i,j) \in \Z^2 \colon 1 \leq i \leq \ell, 1 \leq j \leq \lambda_i}$ be the corresponding \emph{Young diagram}. The \emph{conjugate partition} $\lambda'$ is defined by $\lambda_j' = \max \{i : \lambda_i \geq j\}$. We will freely use implications such as $i \leq j \Rightarrow \lambda_i \geq \lambda_j$. The \emph{hook} $H_{\s z} \subseteq [\lambda]$ is the set of squares weakly to the right and below of $\s z = (i,j) \in [\lambda]$, and the {\it hook length} $h_{\s z} = h_{ij} = |H_{\s z}|= \lambda_i +\lambda'_j - i - j +1$ is the size of the hook. See Figure \ref{fig9}, left drawing.

\medskip

A \emph{standard Young tableau} of shape $\lambda$ is a bijective map $f: [\lambda] \to \{1,\dots,n\}$, such that $f(i_1,j_1) < f(i_2,j_2)$ whenever $i_1 \leq i_2$, $j_1 \leq j_2$, and $(i_1,j_1) \neq (i_2,j_2)$. See Figure \ref{fig9}, right drawing. We denote the number of standard Young tableaux of shape $\lambda$  by $f^\lambda$. The hook-length formula states that if $\lambda$ is a partition of $n$, then
$$f^\lambda = \frac{n!}{\prod_{\s z \in [\lambda]} h_{\s z}}.$$

For example, for $\lambda=(3,2,2)\vdash 7$, the hook-length formula gives
$$f^{322} \, = \, \frac{7!}{5\cdot 4 \cdot 3 \cdot 2\cdot 2 \cdot 1 \cdot 1} \, = \, 21.$$

One way to prove the hook-length formula is by induction on $n$. Namely, it is obvious that in a standard Young tableau, $n$ must be in one of the \emph{corners}, squares $(i,j)$ of $[\lambda]$ satisfying $(i+1,j),(i,j+1) \notin [\lambda]$. Therefore
$$f^\lambda  =  \sum_{\s c \in \p C[\lambda]}  f^{\lambda-\s c},$$
where $\p C[\lambda]$ is the set of all corners of $\lambda$, and $\lambda-\s c$ is the partition whose diagram is $[\lambda] \setminus \set{\s c}$.

\medskip

That means that in order to prove the hook-length formula, we have to prove that $F^\lambda = n!/\prod h_{\s z}$ satisfy the same recursion. It is easy to see that this is equivalent to the following \emph{branching rule for the hook lengths}:
\begin{equation} \label{branch}
 \sum_{ (r,s)  \in  \p C[\lambda]}  \ \, \frac{1}{n}\, \prod_{i=1}^{r-1} \, \frac{h_{is}}{h_{is}-1}\, \prod_{j=1}^{s-1} \, \frac{h_{rj}}{h_{rj}-1}
 \,  =  1.
\end{equation}

In an important development, Green, Nijenhuis and Wilf introduced the \emph{hook walk} which proves \eqref{branch} by a combination of a probabilistic and a short but delicate induction argument \cite{gnw}. Zeilberger converted the hook walk proof into a bijective proof~\cite{Zei}, but laments on the ``enormous size of the input and output'' and ``the recursive nature of the algorithm'' (ibid,~$\S 3$).  With time, several variations of the hook walk have been discovered, most notably the $q$-version of Kerov~\cite{Ker1}, and its further generalization, the $(q,t)$-version of Garsia and Haiman~\cite{GH}. In a recent paper \cite{ckp}, a direct bijective proof of \eqref{branch} is presented. In fact, a bijective proof is presented of the following more general identity, called the \emph{weighted branching formula}:

\begin{equation*} 
\aligned
& \left[\sum_{(p,q) \in [\lambda]} x_p   y_q\right]  \cdot  \left[
\prod_{(i,j) \in [\lambda]\setminus \p C[\lambda]} \, \left(x_{i+1} + \ldots + x_{\lambda_j'}
+  y_{j+1}+ \ldots + y_{\lambda_i}\right)\right] \\
& \, = \, \sum_{(r,s) \in \p C[\lambda]} \left[
\prod_{\stackrel{(i,j) \in [\lambda]\setminus \p C[\lambda]}{\scriptscriptstyle i \neq r,j \neq s}} \, \left(x_{i+1} + \ldots + x_{\lambda_j'}
+  y_{j+1}+ \ldots + y_{\lambda_i}\right)\right]\\
& \, \cdot \,  \left[ \prod_{i=1}^{r} \, \left(x_i + \ldots + x_{r}
+  y_{s+1}+ \ldots + y_{\lambda_i}\right)\right]  \cdot 
\left[ \prod_{j=1}^{s} \, \left(x_{r+1}+ \ldots + x_{\lambda_j'} + y_j + \ldots + y_{s}\right)\right]
\endaligned
\end{equation*}

Here $x_1,\ldots,x_{\ell(\lambda)},y_1,\ldots,y_{\lambda_1}$ are some commutative variables. To see that the stated formula is equivalent to \cite[equation (WHL)]{ckp}, note that in the last products on the right, the terms for $i = r$ and $j = s$ are $x_r$ and $y_s$, respectively. 

\medskip

If we substitute all $x_i$ and $y_j$ by $1$, we get
$$n \cdot \prod_{\s z \in [\lambda] \setminus \p C[\lambda]} (h_{\s z} - 1) \, = \, \sum_{(r,s) \in \p C[\lambda]}\left[\prod_{\stackrel{(i,j) \in [\lambda]\setminus \p C[\lambda]}{\scriptscriptstyle i \neq r,j \neq s}} (h_{\s z} - 1)\right] \prod_{i=1}^r h_{is} \prod_{j=1}^s h_{rj},$$
which is equivalent to \eqref{branch}.

\medskip

Three more identities are also proved, with very similar bijective proofs. Namely, we can replace the sum $\sum_{(p,q) \in [\lambda]} x_p   y_q$ on the left-hand side with $\sum_{p=1}^{\ell(\lambda)} x_p$ and the product $\prod_{j=1}^s$ on the right-hand side with $\prod_{j=2}^s$; we can replace the sum $\sum_{(p,q) \in [\lambda]} x_p y_q$ on the left-hand side with $\sum_{q=1}^{\lambda_1} y_q$ and the product $\prod_{i=1}^r$ on the right-hand side with $\prod_{i=2}^r$; or, we can delete the sum $\sum_{(p,q) \in [\lambda]} x_p   y_q$ on the left-hand side, and replace the products $\prod_{i=1}^r$ and $\prod_{j=1}^s$ on the right-hand side with $\prod_{i=2}^r$ and $\prod_{j=2}^s$, respectively.

\medskip

An open question posed in \cite{ckp} is to find the weighted analogue of the formula
\begin{equation}\label{branchcomp}
 \prod_{\s z \in [\lambda]} (h_{\s z} + 1) \, = \, \sum_{(r,s) \in \p C'[\lambda]}\left[\prod_{\stackrel{(i,j) \in [\lambda]}{\scriptscriptstyle i \neq r,j \neq s}} (h_{\s z} + 1)\right] \prod_{i=1}^{r-1} h_{is} \prod_{j=1}^{s-1} h_{rj}.
\end{equation}
Here $\p C'[\lambda]$ is the set of \emph{outer corners} of $\lambda$, squares $(i,j) \notin [\lambda]$ satisfying $i = 1$ or $(i-1,j) \in [\lambda]$, and $j=1$ or $(i,j-1) \in [\lambda]$. The motivation for this formula is as follows, see \cite{rut}. Division by $\prod_{\s z \in [\lambda]} (h_{\s z} + 1)$ and $\prod_{\s z \in [\lambda]} h_{\s z}$ yields
$$\frac 1{\prod_{\s z \in [\lambda]} h_{\s z}} \, = \, \sum_{(r,s) \in \p C'[\lambda]} \prod_{i=1}^{r-1} \frac 1{h_{is}+1} \prod_{j=1}^{s-1} \frac 1 {h_{rj}+1}\prod_{\stackrel{(i,j) \in [\lambda]}{\scriptscriptstyle i \neq r,j \neq s}} \frac 1{h_{\s z}}$$
We multiply by $(n+1)!$ and use the hook-length formula. We get
$$(n+1) f^{\lambda} = \sum_{\s c \in \p C'[\lambda]} f^{\lambda+\s c},$$
where $\lambda+\s c$ is the partition whose diagram is $[\lambda] \cup \set{\s c}$.

\medskip

Let us introduce the notation $\mu \to \lambda$ or $\lambda \leftarrow \mu$ if $\lambda = \mu - \s c$ for a corner $\s c$ of $\mu$, or, equivalently, if $\mu = \lambda + \s c$ for an outer corner $\s c$ of $\lambda$. We then have
$$\sum_{\mu \vdash n+1} (f^\mu)^2 = \sum_{\mu \vdash n+1} f^\mu \left( \sum_{\lambda \leftarrow \mu} f^\lambda\right) = \sum_{\lambda \vdash n} f^\lambda \left( \sum_{\mu \to \lambda} f^\mu \right) = (n+1) \sum_{\lambda \vdash n} (f^\lambda)^2.$$
Induction proves the famous formula $\sum_{\lambda \vdash n} (f^\lambda)^2 = n!$.

\medskip

This paper is organized as follows. In Section \ref{new}, we present four new formulas. The first is a weighted version of \eqref{branchcomp}, and we call it the \emph{complementary weighted branching rule}. The others are variants of this, similar to the variants of the weighted branching rule presented above. In Section \ref{bijective}, we give a bijective proof of this formula, which is, in particular, the first simple bijective proof of \eqref{branchcomp}. In Section \ref{walks}, we present some results on weighted hook walks, which also give a new way to prove all eight formulas; the proofs of main theorems from this section are deferred to Section \ref{proofs}. And in Section \ref{compl}, we explain how our new formulas can be proved using the four formulas from \cite{ckp} on complementary partitions (which are, roughly, partitions whose diagrams are the complements of $[\lambda]$ in rectangles).

\section{New formulas} \label{new}

The main result of this paper is the following theorem.

\begin{thm}[Complementary weighted branching rule] \label{main}
 Choose a partition $\lambda$, and let $x_1,\ldots,x_{\ell(\lambda)}$, $y_1,\ldots,y_{\lambda_1}$ be some commutative variables. Then
 
 \vspace{-0.4cm}
 
 \begin{equation*} 
\aligned
& \prod_{(i,j) \in [\lambda]} \, \left(x_{i} + \ldots + x_{\lambda_j'}
+  y_{j}+ \ldots + y_{\lambda_i}\right) \, = \!\!\! \sum_{(r,s) \in \p C'[\lambda]} 
\prod_{\stackrel{(i,j) \in [\lambda]}{\scriptscriptstyle i \neq r, j \neq s }} \, \!\!\left(x_{i} + \ldots + x_{\lambda_j'}
+  y_{j}+ \ldots + y_{\lambda_i}\right)\\
& \, \cdot \,  \left[ \prod_{i=1}^{r-1} \, \left(x_{i+1} + \ldots + x_{r-1}
+  y_{s}+ \ldots + y_{\lambda_i}\right)\right]  \cdot 
\left[ \prod_{j=1}^{s-1} \, \left(x_{r}+ \ldots + x_{\lambda_j'} + y_{j+1} + \ldots + y_{s-1}\right)\right].
\endaligned
\end{equation*}
\end{thm}

\begin{figure}[hbt]
\begin{center}
\psfrag{x1}{$x_{1}$}
\psfrag{x2}{$x_{2}$}
\psfrag{x3}{$x_{3}$}
\psfrag{x4}{$x_{4}$}
\psfrag{x5}{$x_{5}$}
\psfrag{y1}{$y_{1}$}
\psfrag{y2}{$y_{2}$}
\psfrag{y3}{$y_{3}$}
\psfrag{y4}{$y_{4}$}
\psfrag{y5}{$y_{5}$}
\psfrag{y6}{$y_{6}$}
\epsfig{file=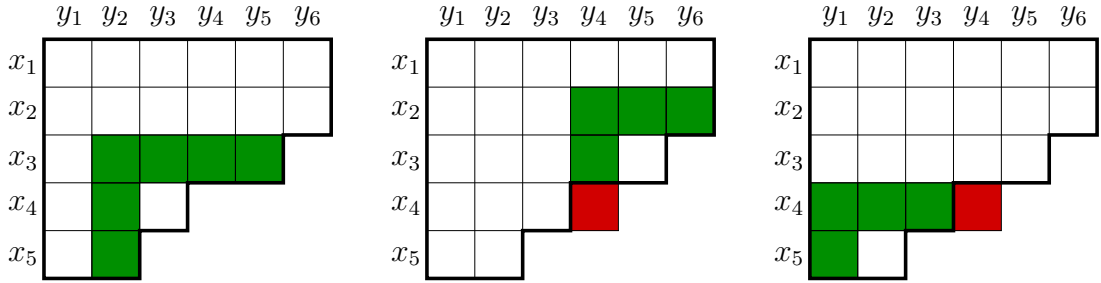,width=0.9 \textwidth}
\end{center}
\caption{The term on the left-hand side of CWBR corresponding to the square $(3,2)$ is $x_3+x_4+x_5+y_2+y_3+y_4+y_5$ (left); the terms on the right-hand side of CWBR corresponding to the outer corner $(4,4)$ and the squares $(2,4)$ and $(4,1)$ are $x_3+y_4+y_5+y_6$ (middle) and $x_4+x_5+y_2+y_3$ (right).}
\label{fig1}
\end{figure}

We refer to this result as CWBR. Figure \ref{fig1} should help understand what the formula is saying. A term on the left-hand side corresponds to a square $(i,j)$ of the diagram $[\lambda]$ and is a natural weighted version of $h_{ij}+1$, see the left diagram in Figure \ref{fig1}. On the right-hand side of CWBR, we have a sum over outer corners. If the square $(i,j)$ of the diagram is in a different row and column of $[\lambda]$ from the chosen outer corner, the corresponding term is the same as on the left-hand side. If it is in the same column, we delete $x_i$, and if it in the same row, we delete $y_j$. Such a term is a weighted version of $h_{ij}$. See the middle and right diagrams of Figure \ref{fig1}. 

\medskip

\begin{exm}
 For $\lambda = 3211$, CWBR gives the following equality:
 \begin{alignat*}{2}
  & \,\, \left(x_1+x_2+x_3+x_4+y_1+y_2+y_3\right)   \left(x_1+x_2+y_2+y_3\right)   \left(x_1+y_3\right)   \left(x_2+x_3+x_4+y_1+y_2\right) \\ 
  & \,\, \left(x_2+y_2\right)   \left(x_3+x_4+y_1\right)   \left(x_4+y_1\right) \\
  =  & \,\,  \left(x_2+x_3+x_4+y_1+y_2\right)   \left(x_2+y_2\right)   \left(x_3+x_4+y_1\right)   \left(x_4+y_1\right)   \left(x_1+x_2+x_3+x_4+y_2+y_3\right) \\
  & \,\,  \left(x_1+x_2+y_3\right)   x_1 \\
  +   & \,\,  \left(x_1+x_2+x_3+x_4+y_1+y_2+y_3\right)   \left(x_1+x_2+y_2+y_3\right)   \left(x_3+x_4+y_1\right)   \left(x_4+y_1\right)   y_3 \\
  &  \,\, \left(x_2+x_3+x_4+y_2\right)   x_2 \\
  +   & \,\,  \left(x_1+x_2+x_3+x_4+y_1+y_2+y_3\right)   \left(x_1+y_3\right)   \left(x_2+x_3+x_4+y_1+y_2\right)   \left(x_4+y_1\right) \\
  &  \,\, \left(x_2+y_2+y_3\right)   y_2   \left(x_3+x_4 \right) \\
  +   &  \,\, \left(x_1+x_2+y_2+y_3\right)   \left(x_1+y_3\right)   \left(x_2+y_2\right)   \left(x_2+x_3+x_4+y_1+y_2+y_3\right)   \left(x_3+x_4+y_1+y_2\right) \\
  &  \,\, \left(x_4+y_1\right)   y_1 
 \end{alignat*}
\end{exm}

We give three more formulas involving outer corners.

$$ {\textstyle \left[{\displaystyle \sum_{p = 1}^{\ell(\lambda)}} {\scriptstyle x_p}\right]  
\prod_{(i,j) \in [\lambda], j \neq 1} \! {\scriptstyle \left(x_{i} + \ldots + x_{\lambda_j'}
+  y_{j}+ \ldots + y_{\lambda_i}\right)}}  = \!\!\!\! \sum_{(r,s) \in \p C'[\lambda], s \neq 1} \!\!\!
{\textstyle \prod_{(i,j) \in [\lambda],i \neq r, j \neq 1,s} \! {\scriptstyle \left(x_{i} + \ldots + x_{\lambda_j'}
+  y_{j}+ \ldots + y_{\lambda_i}\right)}}$$
\begin{equation} \label{x}
\hspace*{5.405cm}  \cdot \,  \left[ \prod_{i=1}^{r-1} \, {\scriptstyle \left(x_{i+1} + \ldots + x_{r-1}
+  y_{s}+ \ldots + y_{\lambda_i}\right)}\right]  \cdot 
\left[ \prod_{j=1}^{s-1} \, {\scriptstyle \left(x_{r}+ \ldots + x_{\lambda_j'} + y_{j+1} + \ldots + y_{s-1}\right)}\right]
\end{equation}

$$ {\textstyle \left[{\displaystyle \sum_{q = 1}^{\lambda_1}} {\scriptstyle y_q}\right]  
\prod_{(i,j) \in [\lambda], i \neq 1} \! {\scriptstyle \left(x_{i} + \ldots + x_{\lambda_j'}
+  y_{j}+ \ldots + y_{\lambda_i}\right)}}  = \!\!\!\! \sum_{(r,s) \in \p C'[\lambda], r \neq 1} \!\!\!
{\textstyle \prod_{(i,j) \in [\lambda],i \neq 1,r, j \neq s} \! {\scriptstyle \left(x_{i} + \ldots + x_{\lambda_j'}
+  y_{j}+ \ldots + y_{\lambda_i}\right)}}$$
\begin{equation} \label{y}
\hspace*{5.405cm}  \cdot \,  \left[ \prod_{i=1}^{r-1} \, {\scriptstyle \left(x_{i+1} + \ldots + x_{r-1}
+  y_{s}+ \ldots + y_{\lambda_i}\right)}\right]  \cdot 
\left[ \prod_{j=1}^{s-1} \, {\scriptstyle \left(x_{r}+ \ldots + x_{\lambda_j'} + y_{j+1} + \ldots + y_{s-1}\right)}\right]
\end{equation}

$$ {\textstyle \left[\!{\displaystyle \sum_{(p,q) \notin [\lambda]}} \!{\scriptstyle x_p   y_q}\right] \! 
\prod_{(i,j) \in [\lambda], i,j \neq 1} \! {\scriptstyle \left(x_{i} + \ldots + x_{\lambda_j'}
+  y_{j}+ \ldots + y_{\lambda_i}\right)}} \! = \!\!\!\!\!\!\!\!\!\! \sum_{(r,s) \in \p C'[\lambda], r,s \neq 1} \!\!\!\!\!\!\!\!\!
{\textstyle \prod_{(i,j) \in [\lambda],i \neq 1,r, j \neq 1,s} \! {\scriptstyle \left(x_{i} + \ldots + x_{\lambda_j'}
+  y_{j}+ \ldots + y_{\lambda_i}\right)}}$$
\begin{equation} \label{xy}
\hspace*{5.405cm}  \cdot \,  \left[ \prod_{i=1}^{r-1} \, {\scriptstyle \left(x_{i+1} + \ldots + x_{r-1}
+  y_{s}+ \ldots + y_{\lambda_i}\right)}\right]  \cdot 
\left[ \prod_{j=1}^{s-1} \, {\scriptstyle \left(x_{r}+ \ldots + x_{\lambda_j'} + y_{j+1} + \ldots + y_{s-1}\right)}\right]
\end{equation}

This last formula requires some clarification: the sum on the left-hand side is over all $(i,j)$ such that $1 \leq i \leq \ell(\lambda)$, $1 \leq j \leq \lambda_1$, $(i,j) \notin [\lambda]$. In other words, we could write
$\left(\sum_{p=1}^{\ell(\lambda)} x_p\right)\cdot\left(\sum_{q=1}^{\lambda_1} y_q\right)-  \sum_{(p,q) \in [\lambda]} x_p   y_q $ instead.

\begin{exm}
 For $\lambda = 3211$, the formulas \eqref{x}, \eqref{y} and \eqref{xy} give

\vspace{-0.4cm}

 \begin{alignat*}{2}
 & \left(x_1+x_2+x_3+x_4\right) \left(x_1+x_2+y_2+y_3\right) \left(x_1+y_3\right) \left(x_2+y_2\right) \qquad \qquad \qquad \qquad \qquad \qquad\\
 = \quad & \left(x_2+y_2\right) \left(x_1+x_2+x_3+x_4+y_2+y_3\right) \left(x_1+x_2+y_3\right) x_1 \\
 + \quad & \left(x_1+y_3\right) \left(x_2+y_2+y_3\right) y_2 \left(x_3+x_4\right) \\
 + \quad & \left(x_1+x_2+y_2+y_3\right) y_3 \left(x_2+x_3+x_4+y_2\right) x_2, \\
\\
& \left(y_1+y_2+y_3\right) \left(x_2+x_3+x_4+y_1+y_2\right) \left(x_2+y_2\right) \left(x_3+x_4+y_1\right) \left(x_4+y_1\right) \\
= \quad & \left(x_3+x_4+y_1\right) \left(x_4+y_1\right) y_3 \left(x_2+x_3+x_4+y_2\right) x_2 \\
+ \quad & \left(x_2+x_3+x_4+y_1+y_2\right) \left(x_4+y_1\right) \left(x_2+y_2+y_3\right) y_2 \left(x_3+x_4\right)\\
+ \quad & \left(x_2+y_2\right) \left(x_2+x_3+x_4+y_1+y_2+y_3\right) \left(x_3+x_4+y_1+y_2\right) \left(x_4+y_1\right) y_1, \\
\\
&  \left(x_3 y_2+x_4 y_2+x_2 y_3+x_3 y_3+x_4 y_3\right) \left(x_2+y_2\right) \\
= \quad & y_3 \left(x_2+x_3+x_4+y_2\right) x_2 \\
+ \quad & \left(x_2+y_2+y_3\right) y_2 \left(x_3+x_4\right).
\end{alignat*}
\end{exm}

\section{Bijective proof of complementary weighted branching rule} \label{bijective}

A direct bijective proof of Theorem \ref{main} shares many characteristics with the bijective proof of the weighted branching rule in \cite[\S 2]{ckp}. We interpret both left-hand and right-hand sides as labelings of the diagram; we start the bijection with a (variant of the) hook walk; and the hook walk determines a relabeling of the diagram. There are, however, some important differences. First, the walk always starts in the square $(1,1)$. Second, the hook walk can never pass through a square that is not in the same row as an outer corner and the same column as an outer corner. Third, the rule for one step of the hook walk is different from the one in \cite{ckp}. And finally, there is an extra shift in the relabeling process. 

\medskip

For the left-hand side of CWBR, we are given a label $x_k$ for some $i \leq k \leq \lambda'_j$, or $y_l$ for some $j \leq l \leq \lambda_i$, for every square $(i,j) \in [\lambda]$. Denote by $F$ the resulting arrangement of $n$ labels (see Figure \ref{fig7}, left), and by $\p F_\lambda$ the set of such labeling arrangements $F$. 

\medskip

For the right-hand side of CWBR, we are given
\begin{itemize}
 \item an outer corner $(r,s)$;
 \item a label $x_k$ for some $i \leq k \leq \lambda'_j$, or $y_l$ for some $j \leq l \leq \lambda_i$, in every square $(i,j) \in [\lambda]$ satisfying $i \neq r$, $j \neq s$;
 \item a label $x_k$ for some $i < k \leq \lambda'_j$, or $y_l$ for some $s \leq l \leq \lambda_i$,
 in every square $(i,s)$;
 \item a label $x_k$ for some $r \leq k \leq \lambda'_j$, or $y_l$ for some $j < l \le \lambda_i$,
 in every square $(r,j)$.
\end{itemize}
Denote by $G$ the resulting arrangement of $n$ labels (see Figure \ref{fig8}), and by $\p G_\lambda$ the set of all such labelings $G$. 

\medskip

Our goal is to give a natural bijection $\varphi: \p F_\lambda \to \p G_\lambda$.

\medskip

We start the bijection by constructing a hook walk. In \cite{ckp}, a label $x_k$ in the square $(i,j)$ meant that we moved to square $(k,j)$, and a label $y_l$ meant that we moved to square $(i,l)$. It should be clear that such a simple rule does not work for CWBR. The first reason is that labels $x_i$ and $y_j$ are also allowed, and this would create a loop. Another important reason is that while the right-hand side of CWBR suggests we should end every hook walk in an outer corner, there are squares of $[\lambda]$ from which an outer corner cannot be reached. In the simplest case of $\lambda = a^b$, we have two outer corners, $(1,a+1)$ and $(b+1,1)$. These two squares can be reached with downward and rightward steps only from the first row and first column of $[\lambda]$. Moreover, we can reach both outer corners only from $(1,1)$.

\medskip

We therefore start the hook walk in $(1,1)$ and move only through squares which are in the same row as an outer corner and in the same column as an outer corner. The rule is as follows. If the current square is $(i,j)$ and the label of $(i,j)$ in $F$ is $x_k$ for $i \leq k \leq \lambda'_j$, move to $(i,\lambda_k+1)$. If the label of $(i,j)$ in $F$ is $y_l$ for $j \leq l \leq \lambda'_j$, move to $(\lambda'_l+1,j)$. Note that in each case, we move to a square which is in the same row as an outer corner and the same column as an outer corner. Moreover, $i \leq k$ implies $\lambda_k \leq \lambda_i$ and $j \leq l$ implies $\lambda'_l \leq \lambda'_j$, so the square we move to is either in $[\lambda]$ or is the outer corner to the right or below $(i,j)$. The process continues until we arrive in an outer corner $(r,s)$, see the right drawing in Figure \ref{fig7}.

\begin{exm} \label{exm}
 Take $\lambda = 988666542$ and the label arrangement drawn in Figure \ref{fig7} on the left.
\begin{figure}[hbt]
\psfrag{x1}{$x_{1}$}
\psfrag{x2}{$x_{2}$}
\psfrag{x3}{$x_{3}$}
\psfrag{x4}{$x_{4}$}
\psfrag{x5}{$x_{5}$}
\psfrag{x6}{$x_{6}$}
\psfrag{x7}{$x_{7}$}
\psfrag{x8}{$x_{8}$}
\psfrag{x9}{$x_{9}$}
\psfrag{y1}{$y_{1}$}
\psfrag{y2}{$y_{2}$}
\psfrag{y3}{$y_{3}$}
\psfrag{y4}{$y_{4}$}
\psfrag{y5}{$y_{5}$}
\psfrag{y6}{$y_{6}$}
\psfrag{y7}{$y_{7}$}
\psfrag{y8}{$y_{8}$}
\psfrag{y9}{$y_{9}$}

\begin{center}
 \epsfig{file=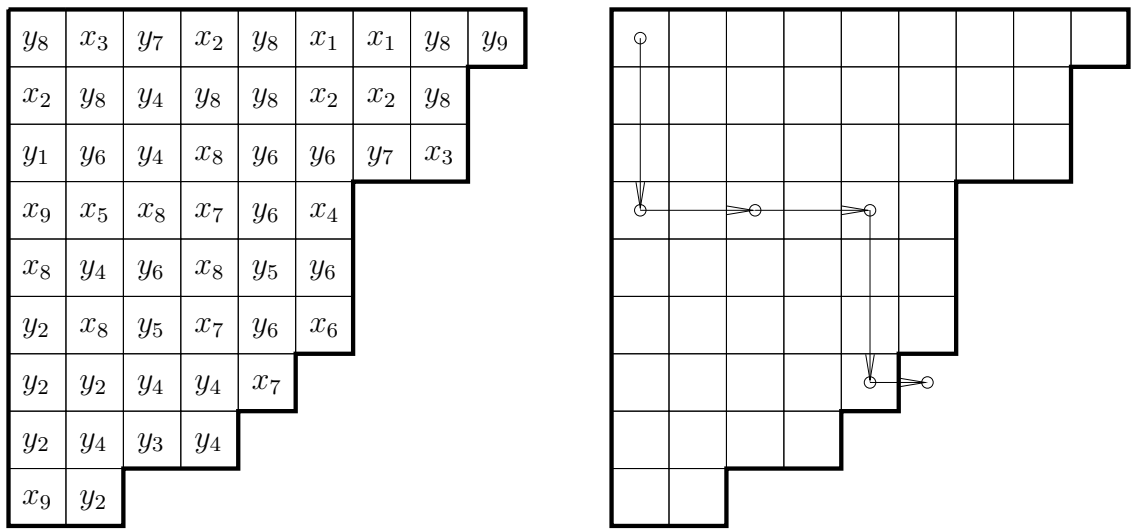,height=7cm}
\end{center}
\caption{An example of an arrangement corresponding to the left-hand side of CWBR for $\lambda = 988666542$; the corresponding hook walk.}
\label{fig7}
\end{figure}

We start in square $(1,1)$. Since the label in $(1,1)$ is $y_8$, we move to $(\lambda'_8+1,1) = (4,1)$. The label in $(4,1)$ is $x_9$, so in the next step, we move to $(4,\lambda_9+1) = (4,3)$. The label there is $x_8$ and our next square is $(4,\lambda_8+1) = (4,5)$. Since the label in $(4,5)$ is $y_6$, we move to $(\lambda'_6+1,5) = (7,5)$. The label in that square is $x_7$ and we therefore move to the outer corner $(7,\lambda_7+1) = (7,6)$ and stop. This hook walk is pictured on the right.
\end{exm}

Shade row $r$ and column $s$. Now we shift the labels in the hook walk and in its projection onto the shaded row and column. If the hook walk has a horizontal step from $(i,j)$ to $(i,j')$, $i \neq r$, move the label in $(i,j)$ right and down to $(r,j')$, and the label from $(r,j)$ up to $(i,j)$. If the hook walk has a vertical step from $(i,j)$ to $(i',j)$, $j \neq s$, move the label from $(i,j)$ down and right to $(i',s)$, and the label from $(i,s)$ left to $(i,j)$. If the hook walk has a horizontal step from $(r,j)$ to $(r,j')$, move the label in $(r,j)$ right to $(r,j')$. If the hook walk has a vertical step from $(i,s)$ to $(i',s)$, move the label in $(i,s)$ down to $(i',s)$. See Figure \ref{fig6}.

\begin{exm}
 We continue with the previous example. On the left, we show how labels trade places. On the right, we have the arrangement after label changes. There are two labels in square $(7,6)$, $x_7$ and $y_6$.
\begin{figure}[hbt]
\psfrag{x1}{$x_{1}$}
\psfrag{x2}{$x_{2}$}
\psfrag{x3}{$x_{3}$}
\psfrag{x4}{$x_{4}$}
\psfrag{x5}{$x_{5}$}
\psfrag{x6}{$x_{6}$}
\psfrag{x7}{$x_{7}$}
\psfrag{x8}{$x_{8}$}
\psfrag{x9}{$x_{9}$}
\psfrag{y1}{$y_{1}$}
\psfrag{y2}{$y_{2}$}
\psfrag{y3}{$y_{3}$}
\psfrag{y4}{$y_{4}$}
\psfrag{y5}{$y_{5}$}
\psfrag{y6}{$y_{6}$}
\psfrag{y7}{$y_{7}$}
\psfrag{y8}{$y_{8}$}
\psfrag{y9}{$y_{9}$}

\begin{center}
 \epsfig{file=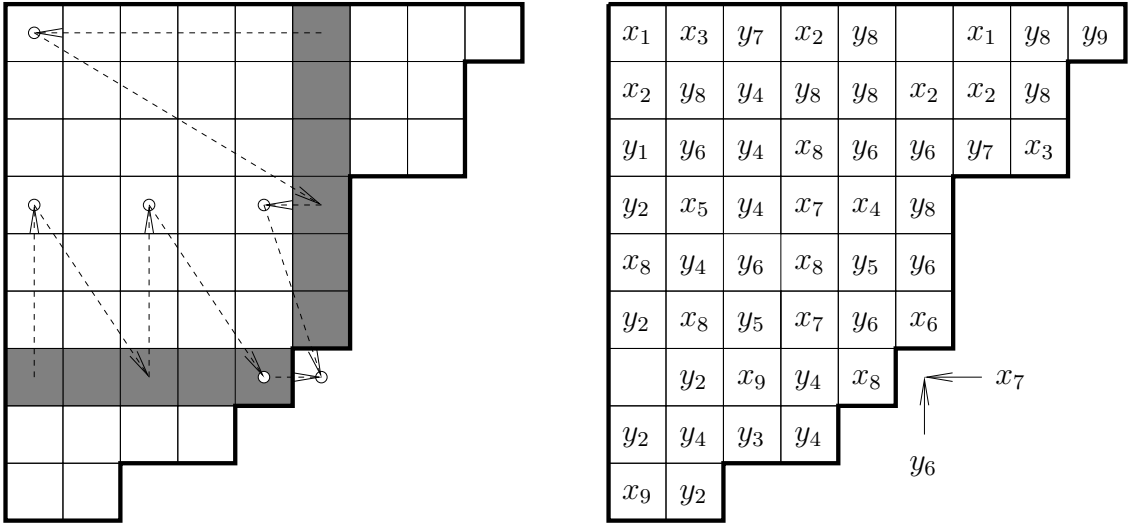,height=7cm}
\end{center}
\caption{An example of the shift of labels for $\lambda = 988666542$.}
\label{fig6}
\end{figure}
\end{exm}

After these changes, we have the following situation. If $r = 1$, there is no label in $(1,1)$, and in $(1,s)$ the label is $x_k$, $1 \leq k \leq \lambda'_{\lambda_1}$. Move all the labels in row $1$ one square to the left. If $s = 1$, there is no label in $(1,1)$, and in $(r,1)$ the label is $y_l$, $1 \leq l \leq \lambda_{\ell(\lambda)}$. Move all the labels in column $1$ one square up. If $r > 1$ and $s > 1$, there are no labels in $(r,1)$ and $(1,s)$. In $(r,s)$, there are two labels: one of the form $x_k$ for $r \leq k \leq \lambda'_{s-1}$, and one of the form $y_l$ for $s \leq l \leq \lambda_{r-1}$. Push all the labels in row $r$, including $x_k$ in $(r,s)$, one square to the left; and push all labels in column $s$, including $y_l$ in $(r,s)$, one square up. See Figure \ref{fig8} for the final arrangement, which we denote $G$.

\begin{exm}
 We continue with the previous example. The following is the final label arrangement.
\begin{figure}[hbt]
\psfrag{x1}{$x_{1}$}
\psfrag{x2}{$x_{2}$}
\psfrag{x3}{$x_{3}$}
\psfrag{x4}{$x_{4}$}
\psfrag{x5}{$x_{5}$}
\psfrag{x6}{$x_{6}$}
\psfrag{x7}{$x_{7}$}
\psfrag{x8}{$x_{8}$}
\psfrag{x9}{$x_{9}$}
\psfrag{y1}{$y_{1}$}
\psfrag{y2}{$y_{2}$}
\psfrag{y3}{$y_{3}$}
\psfrag{y4}{$y_{4}$}
\psfrag{y5}{$y_{5}$}
\psfrag{y6}{$y_{6}$}
\psfrag{y7}{$y_{7}$}
\psfrag{y8}{$y_{8}$}
\psfrag{y9}{$y_{9}$}

\begin{center}
 \epsfig{file=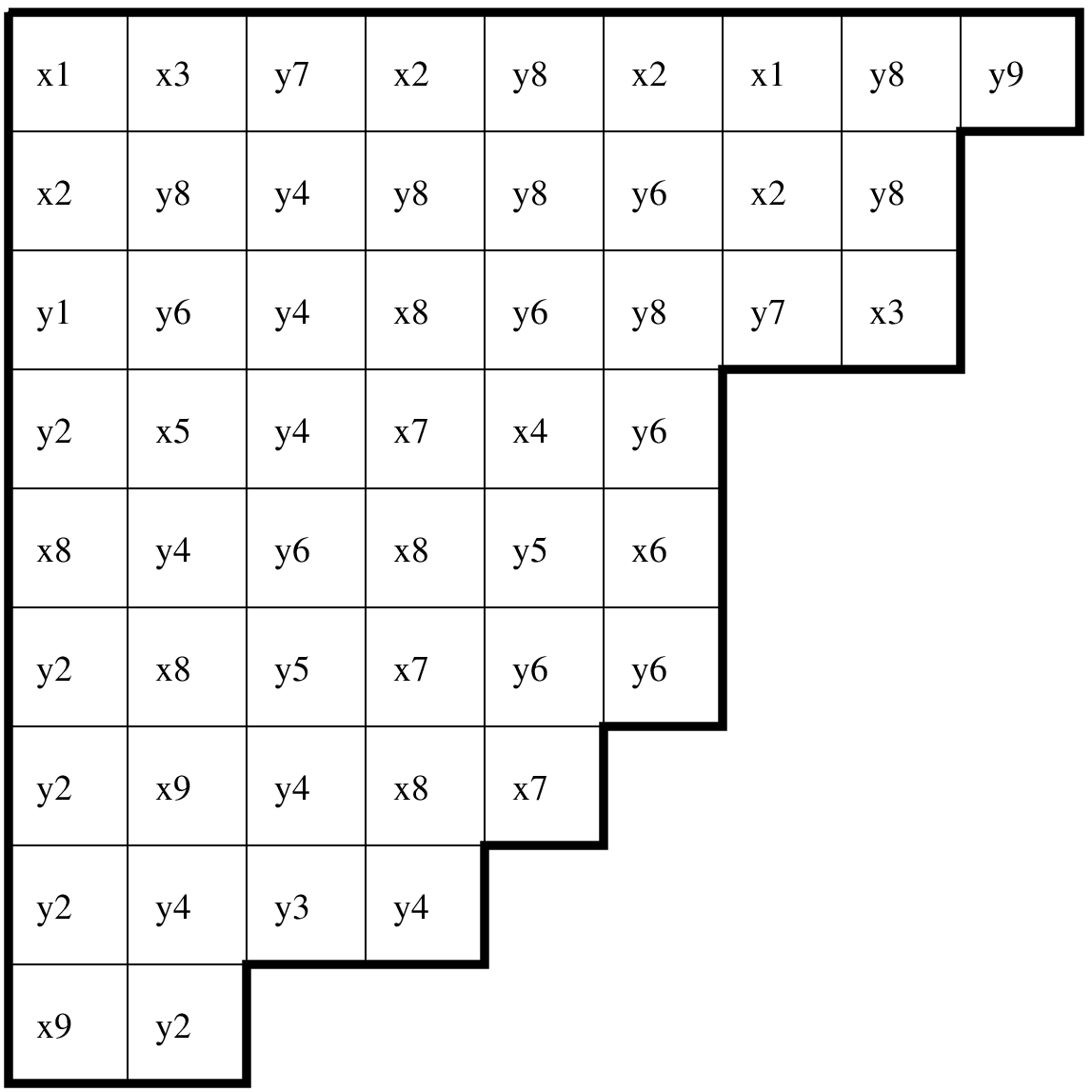,height=7cm}
\end{center}
\caption{The final arrangement.}
\label{fig8}
\end{figure}
\end{exm}

We claim that the final arrangement is in $\p G_\lambda$. If $(i,j)$, $i \neq r$, $j \neq s$, is not one of the squares in the hook walk, then the label of $(i,j)$ in $G$ is the same as in $F$, and it is therefore $x_k$ for $i \leq k \leq \lambda'_j$, or $y_l$ for $j \leq l \leq \lambda_i$. If $(i,j)$ is one of the squares in the hook walk, $i \neq r$, $j \neq s$, then the label of $(i,j)$ in $G$ is the label of either $(i,s)$ or $(r,j)$ in $F$. That means that it is either $x_k$ for $i \leq k \leq \lambda'_s = r-1$, or $y_l$ for $s \leq l \leq \lambda_i$, or $x_k$ for $r \leq k \leq \lambda'_j$, or $y_l$ for $j \leq l \leq \lambda_r = s-1$. In other words, the new label in $(i,j)$ is either $x_k$ for $i \leq k \leq \lambda'_j$, or $y_l$ for $j \leq l \leq \lambda_i$. 

\medskip

The label of $(i,s)$, $1 \leq i \leq r-1$, in $G$ is the label of $(i+1,s)$ in $F$ if $(i+1,s)$ is not in the projection of the hook walk onto column $s$; in other words, it is either $x_k$ for $i +1 \leq k \leq \lambda'_s = r-1$, or $y_l$ for $s \leq l \leq \lambda_{i+1} \leq \lambda_i$. If $(i+1,s)$ is in the projection of the hook walk onto column $s$, we know by construction of the hook walk and relabelings that the new label in $(i,s)$ is $y_l$, where $\lambda'_l = i$. Now $i = \lambda'_l \leq r-1 = \lambda'_s < \lambda'_{s-1}$ implies $l > s - 1$. Also, $l \leq \lambda_{\lambda'_l} = \lambda_i$. In other words, the label in $(i,s)$ is $y_l$ for $s \leq l \leq \lambda_i$. The following is important in the construction of the inverse: since $\lambda'_l = \max \set{k \colon \lambda_k \geq l} = i$, we have $\lambda_{i+1} < l$. In other words, the label in $(i,s)$ is always either $x_k$ for $i +1 \leq k \leq \lambda'_s = r-1$, or $y_l$ for $s \leq l \leq \lambda_i$, and it is $y_l$ for $\lambda_{i+1} < l \leq \lambda_i$ in and only if $(i+1,s)$ is in the projection of the hook walk onto column $s$.

\medskip

We similarly prove that the label in $(r,j)$ is always either $x_k$ for $r \leq k \leq \lambda'_j$, or $y_l$ for $j+1 \leq l \leq \lambda_r = s-1$, and it is $x_k$ for $\lambda'_{j+1} < k \leq \lambda'_j$ in and only if $(r,j+1)$ is in the projection of the hook walk onto row $r$.

\medskip

This shows that $G \in \p G_\lambda$.

\medskip

In the following paragraphs, we sketch the proof of the fact that $\varphi$ has an inverse. The only difficulty lies in reconstructing the hook walk; once we have that, the relabeling process that gives back $F$ is very straightforward.

\medskip

We are given:
\begin{itemize}
 \item an outer corner $(r,s)$;
 \item a label $x_k$ for some $i \leq k \leq \lambda'_j$, or $y_l$ for some $j \leq l \leq \lambda_i$, in every square $(i,j) \in [\lambda]$ satisfying $i \neq r$, $j \neq s$;
 \item a label $x_k$ for some $i < k \leq \lambda'_j$, or $y_l$ for some $s \leq l \leq \lambda_i$,
 in every square $(i,s)$;
 \item a label $x_k$ for some $r \leq k \leq \lambda'_j$, or $y_l$ for some $j < l \le \lambda_i$,
 in every square $(r,j)$.
\end{itemize}

We can read off the projections of the hook walk onto row $r$ immediately. It is the square $(r,1)$, plus all squares $(r,j)$, $j \leq s$, for which the label in $(r,j-1)$ is $x_k$ for $k > \lambda'_j$. Note that since $k \leq  \lambda'_{j-1}$, this can only happen when $(r,j)$ is in the same column as an outer corner. Similarly, the projection of the hook walk onto column $s$ is $(1,s)$ and all squares $(i,s)$, $i \leq r$, for which the label in $(i-1,s)$ is $y_l$ for $l > \lambda_i$.

\medskip

Once we have the projections, it only remains to see whether the hook walk should go right from $(i,j)$, down from $(i,j)$, or terminate. If $i = r$ or $j = s$, the decision is obvious. If $i \neq r$ and $j \neq s$, the label of $(i,j)$ in $G$ is either $x_k$ for $i \leq k \leq \lambda'_j$ or $y_l$ for $j \leq l \leq \lambda_i$. If the label is either $x_k$ for $r \leq k$ or $y_l$ for $l \leq s-1$, we should move to the right; if the label is either $x_k$ for $k \leq r-1$ or $y_l$ for $s \leq l$, move down.

\medskip

We illustrate this with $G$ from the last example, and leave it as an exercise for the reader to check that such a construction indeed gives an inverse of $\varphi$ in general.

\begin{exm}
 Let $G$ be the arrangement in Figure \ref{fig8}, corresponding to the outer corner $(7,6)$. Since the labels of $(7,2)$, $(7,4)$ and $(7,5)$ are $x_9$, $x_8$ and $x_7$, respectively, and since $\lambda'_3 < 9$, $\lambda'_5 < 8$ and $\lambda'_6 < 7$, the projection of the hook walk onto row $7$ contains squares $(7,1)$, $(7,3)$, $(7,5)$ and $(7,6)$. Similarly, since the labels of $(3,6)$ and $(6,6)$ are $y_8$ and $y_6$, respectively, and since $\lambda_4 < 8$ and $\lambda_7 < 6$, the projection of the hook walk onto column $6$ are the squares $(1,6)$, $(4,6)$ and $(7,6)$.\\
 The hook walk starts in $(1,1)$. The label is $x_1$ and $1 \leq 7-1$, so we move down to $(4,1)$. The label there is $y_2$ with $2 \leq 6-1$, so move right to $(4,3)$. The label in $(4,3)$ is $y_4$, and $4 \leq 6-1$. Therefore we move right to $(4,5)$. The label $x_4$ and the inequality $4 \leq 7-1$ imply that we move down to $(7,5)$, and from there we move right to $(7,6)$.\\
 The shifting of labels is easy: move the labels in row $7$ right by one, and the labels in column $6$ down by one. Then reverse the direction of arrows in the right picture in Figure \ref{fig6} and move the labels as indicated by arrows. We get $F$ from Figure \ref{fig7}.
\end{exm}

The proofs of identities \eqref{x}, \eqref{y} and \eqref{xy} are very similar. Note that for an arrangement corresponding to the left-hand side, we now have a chosen row $p$ (respectively, column $q$, respectively, both). We start the hook walk in square $(1,\lambda_p+1)$ (respectively, in $(\lambda'_q+1,1)$, respectively, in $(\lambda'_q+1,
\lambda_p+1)$). It is not difficult to see that such a starting square has second coordinate (respectively, first coordinate, respectively, both coordinates) greater than $1$ and that it is either in $[\lambda]$ or an outer corner. We construct the hook walk in exactly the same fashion as before; we perform the relabeling as before; but before the final shift to the left and up by one, we label $(r,\lambda_p+1)$ (respectively, $(\lambda'_q+1,s)$, respectively, both) with $x_p$ (respectively, $y_q$, respectively, both). The details are left as an exercise for the reader.

\section{Weighted hook walks} \label{walks}

\begin{figure}[hbt]
\begin{center}
\psfrag{R1}{$R_{1}$}
\psfrag{R2}{$R_{2}$}
\psfrag{R3}{$R_{3}$}
\psfrag{R4}{$R_{4}$}
\psfrag{R5}{$R_{5}$}
\psfrag{R6}{$R_{6}$}
\psfrag{R7}{$R_{7}$}
\psfrag{R8}{$R_{8}$}
\psfrag{R9}{$R_{9}$}
\psfrag{R10}{$R_{10}$}
\epsfig{file=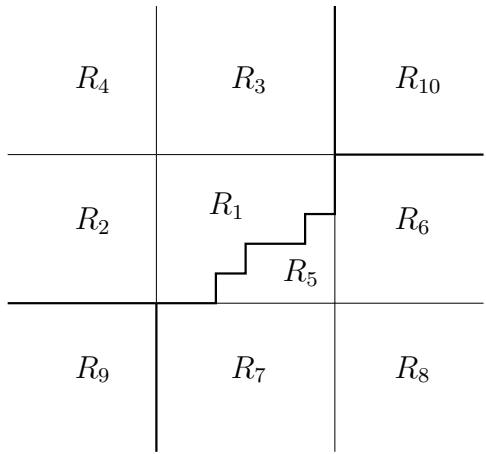,height=6cm}
\end{center}
\caption{Division of the plane into regions $R_1,\ldots,R_{10}$ for $\lambda = 66532$, with some lines in bold.}
\label{fig2}
\end{figure}

Choose a partition $\lambda$ and draw the borders of its diagram in the plane. Now add lines $x = 0$, $x = \ell(\lambda)$, $y = 0$, $y = \lambda_1$; this divides the plane into ten regions $R_1,\ldots,R_{10}$. See Figure \ref{fig2} for an example and the labelings of these regions. Draw the following lines in bold: the half-line $x = 0$, $y \geq \lambda_1$, the half-line $x = \ell(\lambda)$, $y \leq 0$, the half-line $y = 0$, $x \geq \ell(\lambda)$, the half-line $y = \lambda_1$, $x \leq 0$, and the zigzag line separating regions $R_1$ and $R_5$. 

\begin{figure}[hbt]
\begin{center}
\epsfig{file=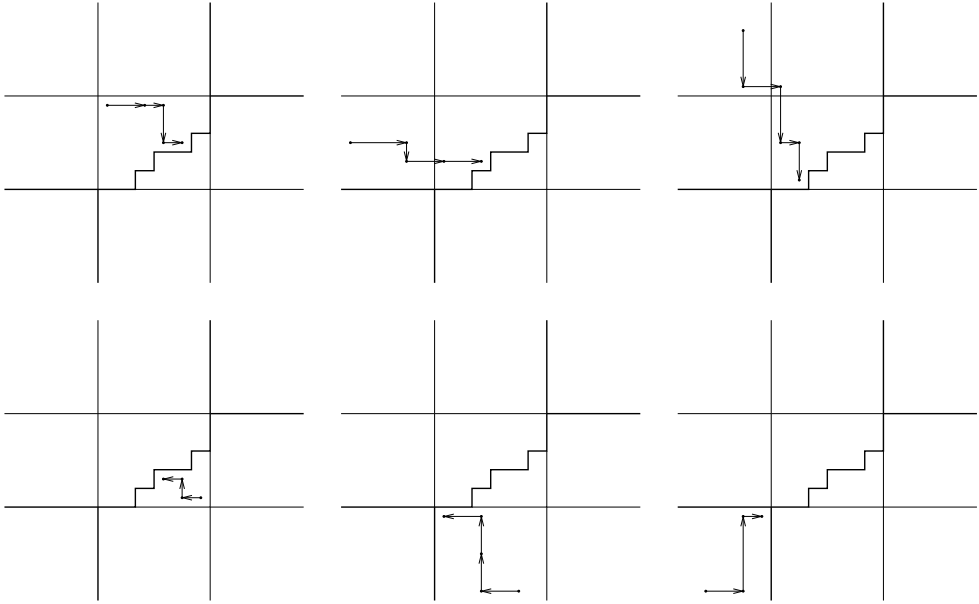,height=8cm}
\end{center}
\caption{Examples of (weighted) hook walks for $\lambda = 66532$.}
\label{fig3}
\end{figure}

Define a \emph{weighted hook walk} as follows. Choose positive weights $(x_i)_{i = -\infty}^{\infty}$, $(y_j)_{j = -\infty}^{\infty}$ satisfying $\sum_i x_i < \infty$, $\sum_j y_j < \infty$. Select the starting square for the hook walk so that the probability of selecting the square $(i,j)$ is proportional to $x_i y_j$. In each step, move in a vertical or horizontal direction toward the bolded line; in regions $R_1$, $R_2,$ $R_3$ and $R_4$, right or down; in regions $R_5$, $R_6$, $R_7$ and $R_8$, left or up; in region $R_9$, right or up; and in region $R_{10}$, left or down. Figure \ref{fig3} shows some examples of weighted hook walks.

\medskip

More specifically, if the current position is $(i,j)$, move to the square $(i',j)$ between $(i,j)$ and the bolded line with probability proportional to $x_{i'}$, and to the square $(i,j')$ between $(i,j)$ and the bolded line with probability proportional to $y_{j'}$. The process stops if we are either in one of the corners of $\lambda$ (if the initial square was in regions $R_1$, $R_2$, $R_3$ or $R_4$), one of the outer corners of $\lambda$ (if the initial square was in regions $R_5$, $R_6$, $R_7$ or $R_8$), the square $(\ell(\lambda)+1,0)$ (if the initial square was in region $R_9$) or $(0,\lambda_1+1)$ (if the initial square was in region $R_{10}$). These last two possibilities are not particularly interesting.

\medskip

The probability of the process ending in a corner $(r,s)$, conditional on starting in $R_1$, was already computed in \cite[Theorem 3]{ckp}. Our goal is to give the probabilities of terminating in a particular corner conditional on starting in $R_2$, $R_3$ and $R_4$, as well as probabilities of ending in a particular outer corner, conditional on starting in $R_5$, $R_6$, $R_7$ and $R_8$. The most interesting observation is that these probabilities turn out to depend only on $x_1,\ldots,x_{\ell(\lambda)}, y_1,\ldots,y_{\lambda_1}$. As a corollary, we obtain the conditional probabilities in the case where all these values are equal. They represent generalizations of classical results due to Greene, Nijenhuis and Wilf from \cite{gnw}, \cite{gnw2}.

\medskip

We extend the definition of $\lambda_i$, $\lambda_j'$ to all $i,j \in \Z$ in a natural way as follows:
\begin{itemize}
 \item for $i \leq 0$, $\lambda_i = \lambda_1$,
 \item for $i \geq \ell(\lambda)+1$, $\lambda_i = 0$,
 \item for $j \leq 0$, $\lambda_j' = \ell(\lambda)$,
 \item for $j \geq \lambda_1 + 1$, $\lambda_j' = 0$.
\end{itemize}

The following statement was proved in \cite{ckp} for $(i,j) \in R_1$ by induction on $|I| + |J|$. 

\begin{lemma} \label{lemma}
 Assume that the weighted hook walk is $(i_1,j_1) \to (i_2,j_2) \to \ldots \to (r,s)$, where $(r,s)$ is either a corner or an outer corner of $\lambda$. Write $I = \{i_1,i_2,\ldots,r\}$ and $J = \{j_1,j_2,\ldots,s\}$ for its vertical and horizontal projections. 
 \begin{enumerate}
  \renewcommand{\labelenumi}{(\alph{enumi})}
  \item Suppose that $(r,s)$ is a corner of $\lambda$. Then the probability that the vertical and horizontal projections are $I$ and $J$, conditional on starting at $(i_1,j_1)$, is
  $$ \textstyle \frac {\prod_{i \in I \setminus \{i_1 \}}x_i} {\prod_{i \in I \setminus \{r \}} (x_{i+1}+\ldots + x_{r}  + y_{s+1} + \ldots + y_{\lambda_i})} \cdot  \frac {\prod_{j \in J \setminus \{j_1\}}y_j} {\prod_{j \in J \setminus \{s\}} (x_{r+1}+\ldots + x_{\lambda'_j} +  y_{j+1} + \ldots + y_{s})}.$$
  \item Suppose that $(r,s)$ is an outer corner of $\lambda$. Then the probability that the vertical and horizontal projections are $I$ and $J$, conditional on starting at $(i_1,j_1)$, is
  $$ \textstyle \frac {\prod_{i \in I \setminus \{i_1 \}}x_i} {\prod_{i \in I \setminus \{r \}} (x_{r}+\ldots + x_{i-1}  + y_{\lambda_i+1} + \ldots + y_{s-1})} \cdot  \frac {\prod_{j \in J \setminus \{j_1\}}y_j} {\prod_{j \in J \setminus \{s\}} (x_{\lambda'_j+1}+\ldots + x_{r-1} +  y_{s} + \ldots + y_{j-1})}.$$
 \end{enumerate}
\end{lemma}
\begin{skt}
 The statement (a) for $(r,s) \in [\lambda]$ is proved in \cite[Lemma 5]{ckp}. If we apply this to the partition whose diagram is $\left(\bigcup_{k=1}^4 R_k\right) \cap \set{(i,j) \colon i \geq \min\set{i_1,1}, j \geq \min\set{j_1,1}}$, we get part (a) in general. Part (b) follows if we rotate the graph by $180^\circ$. \qed
\end{skt}

The following two theorems tell us how to compute probabilities of ending in corners and outer corners. Proofs are deferred to Section \ref{proofs}.

\begin{thm} \label{corners}
 For a corner $\s c = (r,s)$ of $\lambda$, denote by $P(\s c|R)$ the probability that the weighted hook walk terminates in $\s c$, conditional on the starting point being in $R$. Write
 $${\prod}_{rs} = x_r y_s \prod_{i = 1}^{r-1} {\textstyle \left(1+ \frac{x_i}{x_{i+1}+\ldots + x_{r} + y_{s+1} + \ldots + y_{\lambda_i}}\right)} \cdot \prod_{j = 1}^{s-1} {\textstyle \left(1+ \frac{y_j}{x_{r+1}+\ldots + x_{\lambda'_j} + y_{j+1} + \ldots + y_{s}}\right)}.$$
 Then:
 \begin{enumerate}
  \renewcommand{\labelenumi}{(\alph{enumi})}
  \item $P(\s c|R_1)  =  {\frac{1}{\sum_{(p,q) \in [\lambda]}  x_p  y_q}} \cdot {\prod}_{rs}$
  \item $P(\s c|R_2)  =  \frac{1}{\left(\sum_{p=1}^{\ell(\lambda)} x_p\right)(x_{r+1}+\ldots+x_{\ell(\lambda)}+y_1+ \ldots + y_s)} \cdot {\prod}_{rs}$
  \item $P(\s c|R_3) = \frac{1}{\left(\sum_{q=1}^{\lambda_1} y_q\right)(x_{1}+\ldots+x_{r}+y_{s+1}+ \ldots + y_{\lambda_1})} \cdot {\prod}_{rs}$
  \item $P(\s c|R_4) = \frac{1}{(x_{r+1}+\ldots+x_{\ell(\lambda)}+y_1+ \ldots + y_s)(x_{1}+\ldots+x_{r}+y_{s+1}+ \ldots + y_{\lambda_1})} \cdot {\prod}_{rs}$
 \end{enumerate}
 In particular, the sum of each of the above terms over all corners of $\lambda$ equals $1$. Also,
 \begin{enumerate}
  \setcounter{enumi}{4}
  \renewcommand{\labelenumi}{(\alph{enumi})}
  \item $P(\s c) = \frac 1{\left( \sum_{p} x_p\right) \cdot \left( \sum_q y_q\right)}  \cdot\left(1 + \frac{\sum_{p \leq 0} x_p}{x_{1}+\ldots+x_{r}+y_{s+1}+ \ldots + y_{\lambda_1}} \right) \cdot\left(1 + \frac{\sum_{q \leq 0} y_q}{x_{r+1}+\ldots+x_{\ell(\lambda)}+y_1+ \ldots + y_s} \right)  \cdot{\prod}_{rs}$
 \end{enumerate}
\end{thm}

\begin{thm} \label{outercorners}
 For an outer corner, $\s c = (r,s)$ of $\lambda$, denote by $P(\s c|R)$ the probability that the weighted hook walk terminates in $\s c$, conditional on the starting point being in $R$. Write
 $${\prod}'_{rs}= \prod_{i = 1}^{r-1} {\textstyle \left(1 -  \frac{x_i}{x_{i}+\ldots + x_{r-1} + y_{s} + \ldots + y_{\lambda_i}}\right)} \cdot \prod_{j = 1}^{s-1} {\textstyle \left(1- \frac{y_j}{x_{r}+\ldots + x_{\lambda'_j} + y_{j} + \ldots + y_{s-1}}\right)}.$$
 Then:
 \begin{enumerate}
  \renewcommand{\labelenumi}{(\alph{enumi})}
  \item $P(\s c|R_5) = \frac{(x_r + \ldots + x_{\ell(\lambda)} + y_1 + \ldots + y_{s-1})(x_1 + \ldots + x_{r-1} + y_s + \ldots + y_{\lambda_1})}{\sum_{(p,q) \notin [\lambda]} x_p y_q} \cdot {\prod}'_{rs}$
  \item $P(\s c|R_6) = \frac{x_r + \ldots + x_{\ell(\lambda)} + y_1 + \ldots + y_{s-1}}{\sum_{i = 1}^{\ell(\lambda)} x_p} \cdot {\prod}'_{rs}$
  \item $P(\s c|R_7) = \frac{x_1 + \ldots + x_{r-1} + y_s + \ldots + y_{\lambda_1}}{\sum_{q = 1}^{\lambda_1} y_q} \cdot {\prod}'_{rs}$
  \item $P(\s c|R_8) = {\prod}'_{rs}$
 \end{enumerate}
 In particular, the sum of each of the above terms over all outer corners of $\lambda$ equals $1$; note that this proves CWBR, \eqref{x}, \eqref{y} and \eqref{xy}. Also,
 \begin{enumerate}
  \setcounter{enumi}{4}
  \renewcommand{\labelenumi}{(\alph{enumi})}
  \item $P(\s c) = \frac {(x_1 + \ldots + x_{r-1} + \sum_{q=s}^\infty y_q) \cdot (\sum_{p=r}^\infty x_p + y_1 + \ldots + y_{s-1})}{\left( \sum_p x_p\right) \cdot \left( \sum_q y_q\right)} \cdot {\prod}'_{rs}$
 \end{enumerate}
\end{thm}

%
%
%
%

\begin{cor}
 If $x_1 = \ldots = x_{\ell(\lambda)} = y_1 = \ldots = y_{\lambda_1}$, then we have the following. For a corner $\s c = (r,s)$ of $\lambda$,
 \begin{alignat*}{4}
 P(\s c|R_1) & = \frac{f^{\lambda-\s c}}{f^\lambda}, \qquad & P(\s c|R_2) &= \frac{nf^{\lambda-\s c}}{\ell(\lambda)(\ell(\lambda)-r+s)f^\lambda} \\
 P(\s c|R_3) &= \frac{nf^{\lambda-\s c}}{\lambda_1(\lambda_1+r-s)f^\lambda},\qquad & P(\s c|R_4) & = \frac{nf^{\lambda-\s c}}{(\ell(\lambda)-r+s)(\lambda_1+r-s)f^\lambda}
 \end{alignat*}
 In particular, the sum of each of the above terms over all corners of $\lambda$ equals $1$.\\
 For an outer corner, $\s c = (r,s)$ of $\lambda$,
 \begin{alignat*}{4}
 P(\s c|R_5) & = \frac{(\ell(\lambda)-r+s)(\lambda_1+r-s)f^{\lambda+\s c}}{(n+1)(\ell(\lambda)\lambda_1 - n)f^{\lambda}}, \qquad & P(\s c|R_6) &= \frac{(\ell(\lambda)-r+s)f^{\lambda+\s c}}{(n+1)\ell(\lambda) f^{\lambda}}\\
 P(\s c|R_7) &= \frac{(\lambda_1+r-s)f^{\lambda+\s c}}{(n+1)\lambda_1f^{\lambda}},\qquad & P(\s c|R_8) & = \frac{f^{\lambda+\s c}}{(n+1)f^{\lambda}}
 \end{alignat*}
 In particular, the sum of each of the above terms over all outer corners of $\lambda$ equals $1$.
\end{cor}

The corollary gives six new recursive formulas for numbers of standard Young tableaux (and dimensions of irreducible representations of the symmetric group). Recall that one of the classical recursions, $f^\lambda = \sum_{\s c} f^{\lambda-\s c}$, has a trivial bijective proof, and a bijective proof of $(n+1)f^\lambda = \sum_{\s c} f^{\lambda+\s c}$ is essentially the bumping process of the Robinson-Schensted algorithm. It would be nice to find bijective proofs for the new recursions. 

\medskip

Also, we showed in the introduction how the classical recursions prove $\sum_{\lambda \vdash n} \left( f^\lambda\right)^2 = n!$. An interesting question is whether other pairs of ``dual'' recursions, say
$$\ell(\lambda)f^\lambda = n \sum_{\s c} \frac{f^{\lambda-\s c}}{\ell(\lambda)-r+s} \qquad \mbox{and} \qquad (n+1)\ell(\lambda)f^{\lambda} = \sum_{\s c}(\ell(\lambda)-r+s)f^{\lambda+\s c}$$
give a similar identity, and what the version of the Robinson-Schensted proof for that identity would be.

\medskip

The sums over outer corners have the following interesting interpretation. Recall that the content of a square $(i,j)$ of a diagram $[\lambda]$ is defined as $i-j$.

\begin{cor} \label{cor}
 Fix a partition $\lambda \vdash n$. Choose a standard Young tableau of shape $\lambda$ uniformly at random, and an integer $i$, $1 \leq i \leq n+1$ uniformly at random. In the standard Young tableau, increase all integers $\geq i$ by $1$, and use the bumping process of the Robinson-Schensted algorithm to insert $i$ in the tableau. Define the random variable $X$ as the content of the square that is added to $\lambda$. Then
 $$E(X) = 0, \qquad \var(X) = n.$$
\end{cor}
\begin{proof}
 The bumping process is a bijection
 $$\syt(\lambda) \times \set{1,\ldots,n+1} \longrightarrow \bigcup_{\s c \in \p C'[\lambda]} \syt(\lambda+\s c).$$
 This means that the probability that $\s c$ the square added to $\lambda$ is equal to $\frac{f^{\lambda + \s c}}{(n+1)f^\lambda}$. We have
 $$(n+1)\lambda_1 f^{\lambda} = \sum (\lambda_1+r-s)f^{\lambda+\s c} = $$
 $$=\lambda_1 \sum f^{\lambda+\s c} + \sum (r-s)f^{\lambda+\s c} = (n+1)\lambda_1 f^{\lambda} + \sum (r-s)f^{\lambda+\s c}$$
 and therefore
 $$\sum (r-s)f^{\lambda+\s c} = 0,$$
 which is equivalent to $E(X) = 0$. On the other hand, we know that 
 $$ (n+1)(\ell(\lambda)\lambda_1 - n)f^{\lambda} = \sum (\ell(\lambda)-r+s)(\lambda_1+r-s)f^{\lambda+\s c} = $$
 $$=\ell(\lambda) \lambda_1 \sum f^{\lambda+\s c} + (\ell(\lambda)-\lambda_1) \sum (r-s)f^{\lambda+\s c} - \sum (r-s)^2 f^{\lambda+\s c}$$
 and so
 $$\sum (r-s)^2 f^{\lambda+\s c} = (n+1)nf^{\lambda}.$$
 Division by $(n+1)f^{\lambda}$ shows that $\var(X) = n$.
\end{proof}

\begin{rmk}
 The corollary also follows from results of Kerov. From \cite[equations (3.4.3), (3.4.4)]{Ker2}, we get $E(X) = h_1 = p_1$ and $\var(X) = h_2 = p_1^2+p_2/2$, where $p_1=0$ and $p_2 = 2n$ by \cite[equation (3.4.6)]{Ker2}.
\end{rmk}

\section{Proofs of hook walk theorems} \label{proofs}

We only prove parts (d) and (e) of Theorem \ref{corners}, and only part (b) of Theorem \ref{outercorners}, as the proofs of other parts are very similar. 

\medskip

For part (d) of Theorem \ref{corners}, pick $i_1 \leq 0, j_1 \leq 0$, and a corner $\s c = (r,s)$ of $\lambda$. We know that
$$P(\s c|(i_1,j_1)) = \sum_{I,J} P(I,J|(i_1,j_1)),$$
where the sum is over all $I,J$ satisfying $\max I = r$, $\min I = i_1$, $\max J = s$, $\min J = j_1$. By part (a) of Lemma \ref{lemma}, this is
$$\sum_{I,J} \textstyle \frac {\prod_{i \in I \setminus \{i_1 \}}x_i} {\prod_{i \in I \setminus \{r \}} (x_{i+1}+\ldots + x_{r}  + y_{s+1} + \ldots + y_{\lambda_i})} \cdot  \frac {\prod_{j \in J \setminus \{j_1\}}y_j} {\prod_{j \in J \setminus \{s\}} (x_{r+1}+\ldots + x_{\lambda'_j} +  y_{j+1} + \ldots + y_{s})}=$$
$$= \frac{x_r y_s \sum_{I',J'} \frac {\prod_{i \in I'}x_i} {\prod_{i \in I'} (x_{i+1}+\ldots + x_{r}  + y_{s+1} + \ldots + y_{\lambda_i})} \cdot  \frac {\prod_{j \in J'}y_j} {\prod_{j \in J'} (x_{r+1}+\ldots + x_{\lambda'_j} +  y_{j+1} + \ldots + y_{s})}}{(x_{i_1+1} + \ldots + x_r + y_{s+1} + \ldots + y_{\lambda_1})(x_{r+1}+\ldots+x_{\ell(\lambda)}+y_{j_1+1}+\ldots+y_s)},$$
where the sum is over $I' \subseteq \set{i_1+1,\ldots,r-1}$, $J' \subseteq \set{j_1+1,\ldots,s-1}$. It is clear that this is equal to
$$\frac{x_r y_s \prod_{i = i_1+1}^{r-1} \left( 1 + \frac{x_i}{x_{i+1}+\ldots + x_{r}  + y_{s+1} + \ldots + y_{\lambda_i}}\right) \prod_{j = j_1+1}^{s-1} \left(1 + \frac{y_j}{x_{r+1}+\ldots + x_{\lambda'_j} +  y_{j+1} + \ldots + y_{s}}\right)}{(x_{i_1+1} + \ldots + x_r + y_{s+1} + \ldots + y_{\lambda_1})(x_{r+1}+\ldots+x_{\ell(\lambda)}+y_{j_1+1}+\ldots+y_s)}.$$
Now note that
$$\textstyle \frac 1{(x_{i_1+1} + \ldots + x_r + y_{s+1} + \ldots + y_{\lambda_1})} \prod_{i = i_1+1}^{0} \left( 1 + \frac{x_i}{x_{i+1}+\ldots + x_{r}  + y_{s+1} + \ldots + y_{\lambda_i}}\right)=$$
$$\textstyle  =\frac 1{(x_{i_1+1} + \ldots + x_r + y_{s+1} + \ldots + y_{\lambda_1})} \prod_{i = i_1+1}^{0} \frac{x_i+x_{i+1}+\ldots + x_{r}  + y_{s+1} + \ldots + y_{\lambda_1}}{x_{i+1}+\ldots + x_{r}  + y_{s+1} + \ldots + y_{\lambda_1}}=\frac 1{x_{1}+\ldots+x_{r}+y_{s+1}+ \ldots + y_{\lambda_1}}.$$
Together with a similar computation for 
$$\textstyle \frac 1 {x_{r+1}+\ldots+x_{\ell(\lambda)}+y_{j_1+1}+\ldots+y_s} \prod_{j = j_1+1}^{0} \left(1 + \frac{y_j}{x_{r+1}+\ldots + x_{\lambda'_j} +  y_{j+1} + \ldots + y_{s}}\right),$$
this proves that $P(\s c|(i_1,j_1)) = \frac{1}{(x_{r+1}+\ldots+x_{\ell(\lambda)}+y_1+ \ldots + y_s)(x_{1}+\ldots+x_{r}+y_{s+1}+ \ldots + y_{\lambda_1})} \cdot {\prod}_{rs}$. This proves (d).

\medskip

We have
\begin{alignat*}{4}
 P(R_1) & = \textstyle \frac{\sum_{(p,q) \in [\lambda]}  x_p  y_q}{\left( \sum_p x_p\right) \cdot \left( \sum_q y_q\right)}, \qquad & P(R_2) & = \textstyle \frac{\left(\sum_{p=1}^{\ell(\lambda)} x_p\right)\left(\sum_{q \leq 0} y_q\right)}{\left( \sum_p x_p\right) \cdot \left( \sum_q y_q\right)},\\
 \textstyle P(R_3) & = \textstyle \frac{\left(\sum_{p \leq 0} x_p\right)\left(\sum_{q=1}^{\lambda_1} y_q\right)}{\left( \sum_p x_p\right) \cdot \left( \sum_q y_q\right)}, \qquad & P(R_4) & = \textstyle\frac{\left(\sum_{p \leq 0} x_p\right)\left(\sum_{q \leq 0} y_q\right)}{\left( \sum_p x_p\right) \cdot \left( \sum_q y_q\right)},
\end{alignat*}
and therefore, assuming (a)--(d),
\begin{alignat*}{2}
P(\s c) & = P(\s c|R_1) P(R_1) + P(\s c|R_2) P(R_2) + P(\s c|R_3) P(R_3) + P(\s c|R_4) P(R_4) \\
&= \textstyle \frac 1 {\left( \sum_p x_p\right) \cdot \left( \sum_q y_q\right)} \cdot \left(1 + \frac{\sum_{q \leq 0} y_q}{x_{r+1}+\ldots+x_{\ell(\lambda)}+y_1+ \ldots + y_s} + \frac{\sum_{p \leq 0} x_p}{x_{1}+\ldots+x_{r}+y_{s+1}+ \ldots + y_{\lambda_1}}\right. + \\
&+ \textstyle \left. \frac{\left(\sum_{p \leq 0} x_p\right)\left(\sum_{q \leq 0} y_q\right)}{(x_{r+1}+\ldots+x_{\ell(\lambda)}+y_1+ \ldots + y_s)(x_{1}+\ldots+x_{r}+y_{s+1}+ \ldots + y_{\lambda_1})} \right) \cdot {\prod}_{rs},
\end{alignat*}
which is (e).

\medskip

To prove part (b) of Theorem \ref{outercorners}, pick an outer corner $\s c = (r,s)$ of $\lambda$. We want to find
$$P(\s c|R_6) = \frac{\sum_{1 \leq i_1 \leq \ell(\lambda), j_1 > \lambda_1} P(i_1,j_1) \cdot P(\s c|(i_1,j_1))}{P(R_6)}=$$
$$\textstyle =\frac{{\textstyle \sum_{1 \leq i_1 \leq \ell(\lambda), j_1 > \lambda_1}}x_{i_1}y_{j_1} \sum \frac {\prod_{i \in I \setminus \{i_1 \}}x_i} {\prod_{i \in I \setminus \{r \}} (x_{r}+\ldots + x_{i-1}  + y_{\lambda_i+1} + \ldots + y_{s-1})} \cdot  \frac {\prod_{j \in J \setminus \{j_1\}}y_j} {\prod_{j \in J \setminus \{s\}} (x_{\lambda'_j+1}+\ldots + x_{r-1} +  y_{s} + \ldots + y_{j-1})}}{\left(\sum_{p = 1}^{\ell(\lambda)} x_p\right) \left( \sum_{q> \lambda_1} y_q\right)},$$
where the inner sum is over all $I,J$ satisfying $\min I = r$, $\max I = i_1$, $\min J = s$, $\max J = j_1$. We used part (b) of Lemma \ref{lemma} for $P(\s c|(i_1,j_1))$.

\medskip

The trick is to move $x_{i_1}$ into the first inner summation, and to leave $y_{j_1}$ outside. Since
$$\textstyle x_{i_1} \cdot \prod_{i \in I \setminus \{i_1 \}}x_i = x_r \cdot \prod_{i \in I \setminus \{r \}}x_i,$$
we get
$$\textstyle \frac{{\textstyle \sum_{j_1 > \lambda_1}}x_{r}y_{j_1} \left( \sum \frac {\prod_{i \in I \setminus \{r \}}x_i} {\prod_{i \in I \setminus \{r \}} (x_{r}+\ldots + x_{i-1}  + y_{\lambda_i+1} + \ldots + y_{s-1})}\right) \cdot  \left(\sum \frac {\prod_{j \in J \setminus \{j_1\}}y_j} {\prod_{j \in J \setminus \{s\}} (x_{\lambda'_j+1}+\ldots + x_{r-1} +  y_{s} + \ldots + y_{j-1})}\right)}{\left(\sum_{p = 1}^{\ell(\lambda)} x_p\right) \left( \sum_{q > \lambda_1} y_q\right)},$$
where the first inner sum is over all $I$ satisfying $\min I = r$, $\max I \leq \ell(\lambda)$, and the second inner sum is over all $J$ satisfying $\min J = s$, $\max J = j_1$.

\medskip

Let us deal with the inner sums individually. First, we have
$$\sum_{\min I = r,\max I \leq \ell(\lambda)} {\textstyle \frac {\prod_{i \in I \setminus \{r \}}x_i} {\prod_{i \in I \setminus \{r \}} (x_{r}+\ldots + x_{i-1}  + y_{\lambda_i+1} + \ldots + y_{s-1})}} = \prod_{i=r+1}^{\ell(\lambda)} \textstyle \left(1+\frac{x_i}{x_{r}+\ldots + x_{i-1}  + y_{\lambda_i+1} + \ldots + y_{s-1}}\right),$$
and
$$\textstyle x_r \cdot \prod_{i=r+1}^{\ell(\lambda)} \left(1+\frac{x_i}{x_{r}+\ldots + x_{i-1}  + y_{\lambda_i+1} + \ldots + y_{s-1}}\right) = x_r \cdot \prod_{i=r+1}^{\ell(\lambda)} \frac{x_{r}+\ldots + x_{i}  + y_{\lambda_i+1} + \ldots + y_{s-1}}{x_{r}+\ldots + x_{i-1}  + y_{\lambda_i+1} + \ldots + y_{s-1}}=$$
$$= \textstyle x_r \cdot \frac{\prod_{i=r+1}^{\ell(\lambda)} \left(x_{r}+\ldots + x_{i}  + y_{\lambda_i+1} + \ldots + y_{s-1}\right)}{\prod_{i=r+1}^{\ell(\lambda)}\left(x_{r}+\ldots + x_{i-1}  + y_{\lambda_{i}+1} + \ldots + y_{s-1}\right)} = x_r \cdot \frac{\prod_{i=r+1}^{\ell(\lambda)} \left(x_{r}+\ldots + x_{i}  + y_{\lambda_i+1} + \ldots + y_{s-1}\right)}{\prod_{i=r}^{\ell(\lambda)-1}\left(x_{r}+\ldots + x_{i}  + y_{\lambda_{i+1}+1} + \ldots + y_{s-1}\right)} = $$
$$\textstyle = (x_r + \ldots + x_{\ell(\lambda)} + y_1 + \ldots + y_{s-1}) \cdot \frac{\prod_{i=r}^{\ell(\lambda)} \left(x_{r}+\ldots + x_{i}  + y_{\lambda_i+1} + \ldots + y_{s-1}\right)}{\prod_{i=r}^{\ell(\lambda)}\left(x_{r}+\ldots + x_{i}  + y_{\lambda_{i+1}+1} + \ldots + y_{s-1}\right)}=$$
$$=\textstyle (x_r + \ldots + x_{\ell(\lambda)} + y_1 + \ldots + y_{s-1}) \cdot \prod_{i=r}^{\ell(\lambda)} \prod_{j=\lambda_{i+1}+1}^{\lambda_i} \frac{x_r + \ldots + x_i + y_{j+1}+\ldots+y_{s-1}}{x_r + \ldots + x_i + y_{j}+\ldots+y_{s-1}},$$
where the last equality is proved by telescoping. But we have $\lambda_{i+1} < j \leq \lambda_i$ if and only if $i = \lambda'_j$, so reversing the order of multiplication yields
\begin{alignat*}{2}
& (x_r + \ldots + x_{\ell(\lambda)} + y_1 + \ldots + y_{s-1}) \cdot \prod_{j \colon r \leq \lambda'_j \leq \ell(\lambda)} {\textstyle \left( \frac{x_r + \ldots + x_{\lambda'_j} + y_{j+1}+\ldots+y_{s-1}}{x_r + \ldots + x_{\lambda'_j} + y_{j}+\ldots+y_{s-1}}\right)} = \\
=\quad &(x_r + \ldots + x_{\ell(\lambda)} + y_1 + \ldots + y_{s-1}) \cdot \prod_{j=1}^{s-1} {\textstyle \left(1 - \frac{y_j}{x_r + \ldots + x_{\lambda'_j} + y_{j}+\ldots+y_{s-1}} \right)}.
\end{alignat*}
The second computation is very similar. If $s = j_1$, we have
$$\sum_{\min J = s, \max J = j_1} \textstyle \frac {\prod_{j \in J \setminus \{j_1\}}y_j} {\prod_{j \in J \setminus \{s\}} (x_{\lambda'_j+1}+\ldots + x_{r-1} +  y_{s} + \ldots + y_{j-1})} = 1.$$
Otherwise, it is equal to
$${\textstyle \frac{y_s}{x_{\lambda'_{j_1}+1}+\ldots + x_{r-1} +  y_{s} + \ldots + y_{j_1-1}}}\cdot \prod_{j=s+1}^{j_1-1}\textstyle \left(1+\frac {y_j} {x_{\lambda'_j+1}+\ldots + x_{r-1} +  y_{s} + \ldots + y_{j-1}}\right).$$
In either case, we can write this as
$$\textstyle \frac{\prod_{j=s}^{j_1-1} \left(x_{\lambda'_j+1}+\ldots + x_{r-1} +  y_{s} + \ldots + y_{j}\right)}{\prod_{j=s+1}^{j_1} \left(x_{\lambda'_j+1}+\ldots + x_{r-1} +  y_{s} + \ldots + y_{j-1}\right)}=\frac{\prod_{j=s+1}^{j_1} \left(x_{\lambda'_{j-1}+1}+\ldots + x_{r-1} +  y_{s} + \ldots + y_{j-1}\right)}{\prod_{j=s+1}^{j_1} \left(x_{\lambda'_j+1}+\ldots + x_{r-1} +  y_{s} + \ldots + y_{j-1}\right)},$$
and telescoping helps us to write this as
$$\textstyle \prod_{j=s+1}^{j_1} \prod_{i = \lambda'_j+1}^{\lambda'_{j-1}} \frac{x_{i+1}+\ldots + x_{r-1} +  y_{s} + \ldots + y_{j-1}}{x_{i}+\ldots + x_{r-1} +  y_{s} + \ldots + y_{j-1}} = \prod_{i=1}^{r-1} \frac{x_{i+1}+\ldots + x_{r-1} +  y_{s} + \ldots + y_{\lambda_i}}{x_{i}+\ldots + x_{r-1} +  y_{s} + \ldots + y_{\lambda_i}},$$
where we used the fact that $\lambda'_j < i \leq \lambda'_{j-1}$ if and only if $j-1 = \lambda_i$, and that $\lambda'_{j_1}=0$. 

\medskip

Putting these calculations together, our final result for $P(\s c|R_6)$ is
$$\textstyle \frac{{\textstyle \sum_{j_1 > \lambda_1}}\left(y_{j_1} (x_r + \ldots + x_{\ell(\lambda)} + y_1 + \ldots + y_{s-1})\prod_{i=1}^{r-1} \big(1-\frac{x_i}{x_{i}+\ldots + x_{r-1} +  y_{s} + \ldots + y_{\lambda_i}}\big)\prod_{j=1}^{s-1}  \big(1 - \frac{y_j}{x_r + \ldots + x_{\lambda'_j} + y_{j}+\ldots+y_{s-1}}\big)\right)}{\left(\sum_{p = 1}^{\ell(\lambda)} x_p\right) \left( \sum_{q > \lambda_1} y_q\right)}=$$
$$=\frac{x_r + \ldots + x_{\ell(\lambda)} + y_1 + \ldots + y_{s-1}}{\sum_{i = 1}^{\ell(\lambda)} x_p} \cdot {\prod}'_{rs}.$$

\section{Proofs via complementary partitions} \label{compl}

\begin{figure}[hbt]
\begin{center}
\epsfig{file=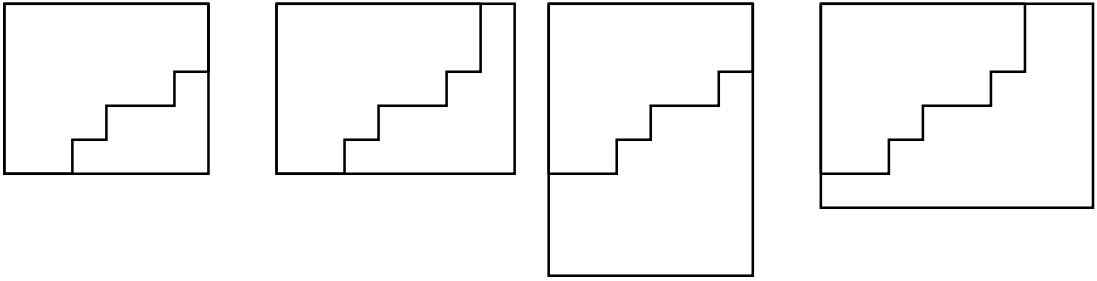,width = 0.9 \textwidth}
\end{center}
\caption{Complementary partitions of $\lambda = 66532$ with respect to $(5,6)$, $(5,7)$, $(8,6)$ and $(6,8)$ are $431$, $54211$, $666431$, $865322$, respectively.}
\label{fig4}
\end{figure}

A partition $\lambda=(\lambda_1,\lambda_2,\ldots,\lambda_\ell)$ of $n$ has several complementary partitions determined by rectangles that contain $[\lambda]$ and have one vertex in $(0,0)$. Namely, choose $a \geq \ell(\lambda)$ and $b \geq \lambda_1$. Pick the non-zero entries of
$$(\underbrace{b,\ldots,b}_{a - \ell},b-\lambda_\ell,b-\lambda_{\ell-1},\ldots b-\lambda_2,b-\lambda_1).$$
We obtain a partition of $ab - n$, which we call the \emph{complementary partition of $\lambda$ with respect to $(a,b)$}. Figure \ref{fig4} represents four different complementary partitions.

It turns out that the formulas CWBR, \eqref{x}, \eqref{y}, \eqref{xy} are equivalent to the four formulas from \cite{ckp} for complementary partitions. We sketch the proof of this statement for \eqref{y} in this section.

\medskip

First note that in \eqref{y}, some terms cancel out. For example, for $\lambda = 3211$, the term $x_4 + y_1$ appears on the left (corresponding to the square $(4,1)$), as well as in all the terms on the right (corresponding to the squares $(4,1)$ for outer corner $(2,3)$, $(4,1)$ for outer corner $(3,2)$, and $(3,1)$ for outer corner $(5,1)$). In general, define $I = \set{i \colon i > 1, (i,s) \in \p C'[\lambda] \mbox{ for some } s}$, $J = \set{j \colon (r,j) \in \p C'[\lambda] \mbox{ for some } r>1}$. Note that $|I| = |J| = |\p C[\lambda]|$. We claim that for $(i,j) \in [\lambda]$, $i>1$, the term $x_{i} + \ldots + x_{\lambda_j'} +  y_{j}+ \ldots + y_{\lambda_i}$ appears (exactly once) in all the terms on the right-hand side of \eqref{y} whenever $i \notin I$ or $j \notin J$.

\medskip

If $i \notin I$ and $j \notin J$, then in particular $i \neq r$ and $j \neq s$ for an outer corner $(r,s)$, so the term $x_{i} + \ldots + x_{\lambda_j'} +  y_{j}+ \ldots + y_{\lambda_i}$ appears in the first product on the right-hand side (and it does not appear in other products, for those, either the lowest $x$-term is $x_r$, or the lowest $y$-term is $y_s$). If $i = r$ and $j \notin J$, then $x_{r} + \ldots + x_{\lambda_j'} +  y_{j}+ \ldots + y_{s-1}$ does not appear in either the first or second product on the right. Since $j \notin J$, we have $\lambda'_{j-1} = \lambda'_j$, and therefore
$$x_{r} + \ldots + x_{\lambda_j'} +  y_{j}+ \ldots + y_{s-1} = x_{r} + \ldots + x_{\lambda_{j-1}'} +  y_{(j-1)+1}+ \ldots + y_{s-1}$$
does appear in the third product on the right. The reasoning for $i \notin I$ and $j = s$ is very similar.

\medskip

This means that \eqref{y} is equivalent to
$$ {\textstyle \left[{\displaystyle \sum_{q = 1}^{\lambda_1}} {\scriptstyle y_q}\right]  
\prod_{(i,j) \in [\lambda] \cap I \times J, i \neq 1} \! {\scriptstyle \left(x_{i} + \ldots + x_{\lambda_j'}
+  y_{j}+ \ldots + y_{\lambda_i}\right)}}  = \!\!\!\! \sum_{(r,s) \in \p C'[\lambda], r \neq 1} \!\!\!
{\textstyle \prod_{\stackrel{(i,j) \in [\lambda] \cap I \times J}{\scriptscriptstyle i \neq 1,r, j \neq s}} \! {\scriptstyle \left(x_{i} + \ldots + x_{\lambda_j'}
+  y_{j}+ \ldots + y_{\lambda_i}\right)}}$$
\begin{equation}\label{y2}
\hspace*{1.505cm}  \cdot \,  \left[ \prod_{i+1 \in I \cap \set{2,\ldots,r}} \, {\scriptstyle \left(x_{i+1} + \ldots + x_{r-1}
+  y_{s}+ \ldots + y_{\lambda_i}\right)}\right]  \cdot 
\left[ \prod_{j + 1 \in J \cap \set{2,\ldots,s}} \, {\scriptstyle \left(x_{r}+ \ldots + x_{\lambda_j'} + y_{j+1} + \ldots + y_{s-1}\right)}\right].
\end{equation}

\medskip

On the other hand, we proved in \cite{ckp} and mentioned in Section \ref{intro} that for every partition $\mu$, we have the equality
\begin{equation*} 
\aligned
& \left[\sum_{q=1}^{\mu_1} y_q\right]  \cdot  \left[
\prod_{(i,j) \in [\mu]\setminus \p C[\mu]} \, \left(x_{i+1} + \ldots + x_{\mu_j'}
+  y_{j+1}+ \ldots + y_{\mu_i}\right)\right] \\
& \, = \, \sum_{(r,s) \in \p C[\mu]} \left[
\prod_{\stackrel{(i,j) \in [\mu]\setminus \p C[\mu]}{\scriptscriptstyle i \neq r,j \neq s}} \, \left(x_{i+1} + \ldots + x_{\mu_j'}
+  y_{j+1}+ \ldots + y_{\mu_i}\right)\right]\\
& \, \cdot \,  \left[ \prod_{i=2}^{r} \, \left(x_i + \ldots + x_{r}
+  y_{s+1}+ \ldots + y_{\mu_i}\right)\right]  \cdot 
\left[ \prod_{j=1}^{s} \, \left(x_{r+1}+ \ldots + x_{\mu_j'} + y_j + \ldots + y_{s}\right)\right]
\endaligned
\end{equation*}

Define $I' = \set{i \colon (i,s) \in \p C[\mu] \mbox{ for some } s}$ and $J' = \set{j \colon (r,j) \in \p C[\mu] \mbox{ for some } r}$. We can prove now that a term $x_{i+1} + \ldots + x_{\mu_j'}+  y_{j+1}+ \ldots + y_{\mu_i}$ cancels out from the above equality whenever $i \notin I'$ or $j \notin J'$. That means that we have
\begin{equation*} 
\aligned
& \left[\sum_{q=1}^{\mu_1} y_q\right]  \cdot  \left[
\prod_{(i,j) \in [\mu] \setminus \p C[\mu] \cap I' \times J'} \, \left(x_{i+1} + \ldots + x_{\mu_j'}
+  y_{j+1}+ \ldots + y_{\mu_i}\right)\right] \\
& \, = \, \sum_{(r,s) \in \p C[\mu]} \left[
\prod_{\stackrel{(i,j) \in [\mu]\setminus \p C[\mu] \cap I' \times J'}{\scriptscriptstyle i \neq r,j \neq s}} \, \left(x_{i+1} + \ldots + x_{\mu_j'}
+  y_{j+1}+ \ldots + y_{\mu_i}\right)\right]\\
& \, \cdot \,  \left[ \prod_{\stackrel{\scriptstyle i-1 \in I' \cap}{ \set{1,\ldots,r-1}}} \!\!\!\! \left(x_i + \ldots + x_{r}
+  y_{s+1}+ \ldots + y_{\mu_i}\right)\right]  \cdot 
\left[ \prod_{\stackrel{\scriptstyle j-1 \in J'\cap  }{\set{0,\ldots,s-1}}} \!\!\!\! \left(x_{r+1}+ \ldots + x_{\mu_j'} + y_j + \ldots + y_{s}\right)\right]
\endaligned
\end{equation*}

\medskip

It turns out that if we write this identity for $\mu$ the complement of $\lambda$ with respect to $(\ell(\lambda)+1,\lambda_1)$, with $x_i$ replaced by $x_{\ell(\lambda)+2-i}$, and with $y_j$ replaced by $y_{\lambda_1+1-j}$, we get \eqref{y2}. 

\medskip

\begin{figure}[hbt]
\begin{center}
\epsfig{file=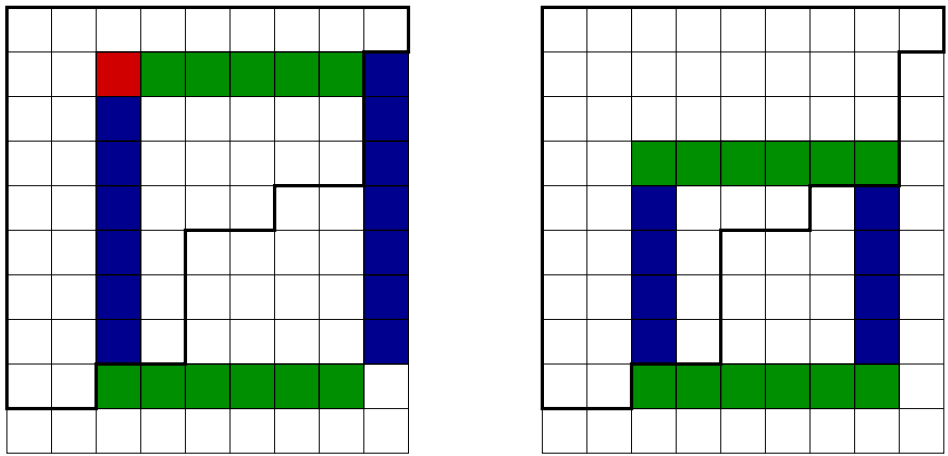,height=6cm}
\end{center}
\caption{Hooks for the partition and its complement when $\lambda=988864442$.}
\label{fig5}
\end{figure}

The geometric reason for that is as follows. If $(i,j)$ is a square of $\lambda$ that is in the same row as an outer corner and the same column as an outer corner, the hook of $(i,j)$ in $[\lambda]$, with $(i,j)$ counted twice, is the same as the hook of $(\lambda'_j+1,\lambda_i+1)$ without the square $(\lambda'_j+1,\lambda_i+1)$ in the complement of $\lambda$ with respect to $(\ell(\lambda)+1,\lambda_1)$, see Figure \ref{fig5}, left. Moreover, $(\lambda'_j+1,\lambda_i+1)$ is in the same row and column as a corner of the complement of $\lambda$ with respect to $(\ell(\lambda)+1,\lambda_1)$. Similarly, the hook of $(i,s)$, $i+1 \in I$, in $\lambda$ is the same as the hook of the square $(r,\lambda_i)$ in the complement of $\lambda$ with respect to $(\ell(\lambda)+1,\lambda_1)$, see Figure \ref{fig5}, right. Moreover, $(r,\lambda_i)$ is next to a square that is in the same column as a corner of the complement of $\lambda$ with respect to $(\ell(\lambda)+1,\lambda_1)$. This is also the reason why the telescoping argument in the previous section worked. We omit the details.

\medskip

For CWBR, we would take the complement of $\lambda$ with respect to $(\ell(\lambda)+1,\lambda_1+1)$; for \eqref{x}, with respect to $(\ell(\lambda),\lambda_1+1)$ and for \eqref{xy}, with respect to $(\ell(\lambda),\lambda_1)$. The details are left to the reader.

\bigskip

{\bf Acknowledgments.} The author would like to thank Paul Edelman and Igor Pak for many helpful comments, and Dan Romik for pointing out how to relate Corollary \ref{cor} to Kerov's work.

\end{document}